\documentclass[11pt]{article}
\usepackage{amsfonts, amsbsy, amssymb, eucal}
\usepackage{graphicx}

\begin{document}

%%%%%%%%%%%%%%%%%%%

\def\L{{\mathbb L}}
\def\B{{\mathbb B}}
\def\R{{\mathbb R}}
\def\Rn{{\mathbb R}^n}
\def\Rnp{{\mathbb R}^{n+1}}
\def\T{{\mathbb T}}
\def\N{{\mathbb N}}
\def\M{{\mathbb M}}
\def\C{{\mathbb C}}
\def\Z{{\mathbb Z}}
\def\Q{{\mathbb Q}}
\newcommand{\ld}{,\ldots,}
\newcommand{\spr}[2]{\langle#1,#2\rangle}
\newcommand{\avr}[1]{\langle#1\rangle}

\newtheorem{definition}{Definition}
\newtheorem{assumption}{Assumption}
\newtheorem{theorem}{Theorem}
\newtheorem{lemma}{Lemma}[theorem]
\newtheorem{corollary}{Corollary}[theorem]

%%%%%%%%%%%%%%%%%%%
\title{Integrability versus topology of configuration manifolds and domains of possible
motions}
\author{M. Rudnev\thanks{Contact address: Department of Mathematics,
University of Bristol, University Walk, Bristol BS8 1TW, UK;\newline
e-mail: {\tt m.rudnev@bris.ac.uk}} $\,$and V. Ten\thanks{Contact
address: Department of Mathematics, University of Bristol,
University Walk, Bristol BS8 1TW, UK;\newline  e-mail: {\tt
v.ten@bris.ac.uk} }}
\date{January 23, 2004}
\maketitle

\begin{abstract}
We establish a generic sufficient condition for a compact $n$-dimensional
manifold bearing an integrable geodesic flow to be the $n$-torus. As a complementary
result, we show that in the case of domains of possible motions with
boundary, the first Betti number of the domain of possible motions
may be arbitrarily large.
\end{abstract}

\def\proclaim#1{{\bf #1.}\bgroup\it}
\def\endproclaim{\egroup}
\def\leqslant{\leq}
\def\Bbb#1{{\bf #1}}

\medskip
\medskip
\noindent {AMS subject classification 37J30, 37J35}

\subsection*{Introduction}

\indent Analytic integrability of natural Hamiltonian systems
implies stringent restrictions on topology of the configuration
manifold. Topological obstructions to integrability were first
introduced by V.V. Kozlov in \cite{Ko}, where he showed that a
compact two-dimensional orientable manifold can be a configuration
space of an integrable natural system only if this manifold is
diffeomorphic to a sphere or a torus.

This result was generalized to the $n$-dimensional case by I.
Taimanov (\cite{Ta}, see also \cite{BT}). The following necessary
conditions for analytic integrability of a geodesic flow on a
$n$-dimensional compact manifold $M^n$ were found.
\begin{enumerate}
\item The fundamental group $\pi_1(M^n)$ is almost commutative. I.e. $\pi_1(M^n)$ contains a commutative
subgroup of a finite index.
\item The first Betti number
$\beta_1(M^n)\leq n$, and the real cohomology ring $H^*(M^n,\R)$
contains a subring, isomorphic to $H^*(\T^k,\R),$ where
$k=\beta_1(M^n)$, and $\T^k$ denotes the $k$-dimensional torus.
\end{enumerate}

Recall that the dynamics for any natural Hamiltonian system gets
reduced to the geodesic flow dynamics via the Jacobi metric, provided that total energy of the system exceeds
the absolute maximum of the potential, i.e. the domain of possible
motions is the whole configuration manifold (and in particular is without boundary).

Observe that Taimanov's theorem implies that if
$\beta_1(M^n)=n$, then the ring $H^*(M^n,\R)$ is isomorphic to that
of the $n$-torus. The theorem in this note establishes the fact that under
a generic assumption of non-degeneracy of the momentum map (precise definitions are given further in
this note), the manifold $M^n$ is actually {\em
diffeomorphic} to $\T^n$, in the case when $\beta_1(M^n)=n.$

The situation turns out to be dramatically different however in the case
when the domain of possible motions does not coincide with the whole
configuration manifold $M^n$, having a non-empty boundary instead.
(For a natural system Hamiltonian this is the case when total
energy is less than the absolute maximum of the potential on $M^n$.)
For bounded domains of possible motions, in the case $n=2$, Kozlov and Ten (\cite{KT}) showed that there
are no obstructions to topology of domains of possible motions with
boundary, for integrable systems. The second result of this note is
a multidimensional generalization of this fact: there is no
restriction on the first Betti number of domains of possible
motions.

\subsection*{Definitions and results}
Let us briefly recall the standard definitions, see
e.g. \cite{A1}.

Let $(M^n,g)$ be a compact
Riemannian manifold with a  non-degenerate metric $g=(g_{ij}(x))$. It defines the geodesic flow on $TM^n$, which
satisfies the Lagrange equations, with the Lagrangian
$$
L(x,\dot x)={1\over2}g_{ij}(x)\dot{x}^i\dot{x}^j.
$$
The Legendre transform $TM^n \mapsto T^*M^n,$ defined via
$$\dot{x}\rightarrow y:\,y_i=g_{ij}(x)\dot{x}^j,$$
results in a Hamiltonian system on $T^*M^n$, with the Hamiltonian
$$
H_g(x,y)={1\over2}g^{ij}(x)y_i y_j, \;\mbox{ where }\;g^{ij}g_{jk}=\delta^{i}_k,
$$
preserving the symplectic form $\omega= dy\wedge dx$. The latter defines the
Poisson bracket for functions on $T^*M^n$:
$$\{f,g\}=\omega(df, dg)={\partial_y f}{\partial_x g} - {\partial_y g}{\partial_x
f}.$$

A natural Hamiltonian system is a system with Hamiltonian
$$
H_\nu (x,y)={1\over2}g^{ij}(x)y_i y_j + U(x),
$$
where $U$ is a potential. Given a value of total energy $E$, exceeding $\sup_M
U$, according to the Maupertuis principle, the phase trajectories of the motions with total energy  $E$
are geodesic lines of the Jacobi metric $g_J= (E-U)g$. Thus the case $E>\sup_M
U$ falls into the realm of geodesic flows. The Jacobi metric becomes
degenerate if $E\leq\sup_M U$, the condition $U(x)\leq E$ defining the domain
of possible motions, with boundary, on the manifold $M^n$.

It is said that a Hamiltonian system on a symplectic manifold $P^{2n}$ (in our case $P^{2n}=T^*M^n$) is (Liouville) integrable if
there are $n$ functionally independent first integrals $F_1,\ldots,F_n$ in
involution, i.e. $\{F_i,F_j\}=0,\,\forall i,j=1,\ldots,n.$ One can always take
$F_1=H$, the Hamiltonian.

\medskip \noindent {\em Remark.} From this point on, we would like to emphasize
that in this note we deal with analytic integrability only,
namely when the objects $M,g,U$ are analytic, in particular when $g_{ij}(x)$
and $U(x)$ are real analytic functions of $x\in M^n$. $C^\infty$ integrability
is strikingly different. E.g. there exist $C^\infty$, but not analytically
integrable flows with positive topological entropy, see \cite{BT}.

\medskip
The map $F:P^{2n}\to\R^n$, where $F(x)=(F_1(x),\dots,F_n(x))$ is
called the momentum map. By the Liouville theorem, under additional
assumption that the differential $dF$ has a full rank $n$ everywhere
on a level set, or {\em layer}, $F_c=\{F(x)=c\}$  in $P^{2n}$, the
layer in question is an $n$-torus, with the linear flow thereupon.

The level set $F_c$ is critical if for some $x\in F_c$, $dF(x)=0$.
Dynamics within critical layers can be quite involved. We will take
advantage of the following topological features regarding the
critical values of the momentum map, see \cite{Fo}, \cite{Fom} for
thorough discussion and further references.
\begin{enumerate}
\item A critical layer $F_c$ is a union of non-intersecting strata.
Each stratum is a smooth invariant manifold.
The closure of a stratum includes the stratum itself
as well as possibly some strata of lower dimensions.
\item
Liouville tori together with stratified critical layers give rise to an
$n$-dimensional cell complex $C$ in the space of momenta $F$, whose
$n$-dimensional cells consist of points, which are images (under the momentum
map) of Liouville tori. Cells of dimension $k<n$
consist of critical values $c=(c_1,\ldots, c_n)$.
\end{enumerate}

Our main assumption for the forthcoming Theorem 2 is non-degeneracy
of the momentum map, as follows.

\renewcommand{\theenumi}{\roman{enumi}}
\medskip\noindent
\proclaim{Definition 1} A system is called {\it non-degenerate} if
\begin{enumerate} \item  the cell complex $C$ is
finite;
\item in every stratum of dimension $k\in\{n,\ldots,1\}$, the rank of the momentum map equals $k$;
\item layers in $P^{2n}$, which are pre-images of points $c$, belonging to a single
cell of the complex $C$ are homotopically equivalent, with
stratification.
\end{enumerate}
\endproclaim

\medskip
\noindent {\em Remark.} Observe that the clause (i) of Definition 1
is fulfilled by analyticity.

\medskip
Let us now state the main results of this note.

\medskip\noindent
\proclaim{Theorem 2} Let $M^n$ be an analytic Riemannian manifold,
with $\beta_1(M^n)=n$. If a geodesic flow on $M^n$ is integrable and
non-degenerate, the manifold $M^n$ is diffeomorphic to
$\T^n$.
\endproclaim

\medskip
The proof of Theorem 2 is given in the next section. Now, in
contrast to the statement of the theorem, let us address the case of
natural systems with domains of possible motions with boundary. We
further describe a construction, which shows that the domain of
possible motions with boundary can possess an arbitrarily large
first Betti number.

\medskip
\noindent
{\bf Example 3.} Consider an $n$-dimensional cube $Q=[-\pi,\pi]^n$.
For $x=(x_1,\ldots,x_n)\in Q$, define a $2\pi$-periodic potential
$$
U(x) =\sum_{i=1}^n\cos{x_i}.
$$
A natural system with Hamiltonian
$$H={1\over2}\sum_{i=1}^n y_i^2 + U(x)$$ is obviously integrable.

Consider the set $\{U\leq n-2\}\subset Q.$ This set intersects each
$n-1$ dimensional face of the cube $Q$ at one point  only, which is
the center of the corresponding face. Let us fix the total energy
value $E=n-2-\epsilon$ and look at the domain of possible motions.
Clearly, the domain would look like an inflated version of the union
of coordinate axes, i.e is homeomorphic to the union of
$\epsilon$-tubular neighborhoods of the coordinate axes, within $Q$.
In fact, it is safe to set $\epsilon=1$, so let us define a manifold
with a boundary $D=\{V\leq n-3\}\subset Q$. Let us refer to the
intersections of $D$ with the faces of $Q$ as faces of $D$. Due to
the $2\pi$-periodicity of the potential, an arbitrary large number
$N$ of building blocks (e.g. $N=2^n$, to fill in $[-\pi,3\pi]^n$ or
$N=3^n$, to fill in $[-3\pi,3\pi]^n$) which are identical copies of
$D$, can be glued along opposite faces. Let us called the resulting
building $B$. After the building has been completed, its opposite
faces are identified. E.g. in the above illustration with $N=2^n$,
the Hamiltonian $H$ should be viewed as a function on $T^*\T_2^n$,
where we denote $\T_2^n\equiv [-\pi,3\pi]^n$. Alternatively, on the
standard torus $\T^n$, one could take the potential as
$\sum_{i=1}^n\cos k x_i,\,k=2,3,\ldots$, or more generally as
$\sum_{i=1}^n\cos k_i x_i$, with $k_i\in \N\equiv\{1,2,\ldots\}.$

The above described procedure clearly enables one to build domains
of possible motion that possess an arbitrarily large first Betti
number.

\subsection*{Proof of Theorem 2}

Every free homotopy class of $M^n$ has a closed geodesic of minimal length.
Hence, let $Z_\alpha$ be the set of all closed geodesic of minimal length of class $\alpha$,
and $L_\alpha$ be their length. Theorem 2 will follow from the following lemmata.

\medskip\noindent
\proclaim{Lemma 4} A stratum of dimension $k\leq n$ is diffeomorphic
to a cylinder $\T^{k-m}\times\R^{m}$, for some $0\leq m < k$.
\endproclaim

\medskip\noindent
\proclaim{Definition 5} Liouville tori such that the natural
projection on $M^n$ of any non-trivial cycle thereupon is
non-trivial are called non-trivial. And trivial otherwise.
\endproclaim

\medskip\noindent
\proclaim{Lemma 6} For some class $\alpha$ of infinite order, the
phase-space pre-images (under the inverse of the natural projection
of $T^*M^n$ onto $M^n$) of all geodesics from the set $Z_{k\alpha},$
belong to {\it non-trivial} Liouville tori, for all
$k\in\N$.\endproclaim

\medskip\noindent
\proclaim{Lemma 7} For any $\alpha,$ every $n$-dimensional cell of
the complex $C,$ containing a non-trivial torus, contains geodesics
of at most one set $Z_{k\alpha},$ $k\in\N$.
\endproclaim

\medskip
Let us show that the lemmata suffices to prove Theorem 2. Let the
class $\alpha$ be provided by Lemma 6. Assuming Lemma 7 implies that
the ratio $d_k={L_{k\alpha}\over k}$ can take only a finite number
of values (as the number of cells in the complex $C$ has been
assumed finite), so $d_k$ reaches its minimum for some $k=k_\alpha$
and the corresponding class $k_\alpha\alpha$. Let $T$ be a
non-trivial (by Lemma 6) torus carrying the phase space pre-image of
a geodesic from the set $Z_{k_\alpha\alpha}$. As the flow on
Liouville tori is linear, all the orbits on the torus $T$ have equal
length, hence belong to the same set $Z_{k_\alpha\alpha}$.

It is easy to see that the projections of distinct orbits on $T$
(alias closed geodesics on $M^n$) either have no common points or
coincide as sets. Indeed, suppose there is a pair of such
non-coinciding intersecting geodesics. Then their homotopy sum has
length $2L_{k_\alpha\alpha}$ and belongs to the class of
$2k_\alpha\alpha$. Therefore
$2L_{k_\alpha\alpha}/2k_\alpha=d_{k_\alpha}$. It follows that the
homotopy sum in question has minimum possible length. This is a
contradiction, as on the other hand, its length can be decreased, as
non-coinciding geodesics always intersect transversely.

So, the projections of a pair of orbits on $T$, as subsets of $M^n$,
either coincide or are disjoint as sets. The former side of the
alternative may occur only if some closed geodesic is covered twice,
by a pair of orbits $\gamma=\{(x_\gamma,y_\gamma)\}\subset T$ and
$\gamma'=\{(x_\gamma,-y_\gamma)\}\subset T.$

Consider a small neighborhood $B_{\xi_0}\subset T$ of a point
$\xi_0=(x,y) \in T$. This neighborhood could be taken small enough
to ensure $(x,-y) \not\in B_{\xi_0}$. Hence, the orbits passing
through $B_{\xi_0}$ locally project on $M^n$ one-to-one. Therefore
the configuration space is regularly covered by the $n$-torus $T$.
On the other hand, the degree of covering cannot exceed $2,$ (as it
was shown that the projections of any two orbits from $T$ cannot
intersect transversally in $M^n$, by minimality). Suppose finally,
some closed geodesic is covered twice, by a pair of orbits $\gamma$
and $\gamma'$ as above. This implies that the class $k_\alpha\alpha$
of the geodesic in question is of finite order, namely $2$. (Indeed,
the sum of the classes of the projections of the orbits $\gamma$ and
$\gamma'$ onto $M^n$ is homotopic to zero.) This contradicts the
assumption that the class $\alpha$ is of infinite order. Hence, $T$
covers $M^n$ globally one-to-one, which proves the theorem,
conditional on the lemmata.

\medskip
\noindent {\bf Proof of Lemma 4.} The proof follows after going
through standard the proof of the Liouville-Arnold theorem
(\cite{A1}), using the clause (ii) of Definition 1.

\medskip
\noindent {\bf Proof of Lemma 6.} Assume the contrary, i.e that for
every class $\alpha$ of infinite order, there exists some
$k^*_\alpha,$ such that the set $Z_{k^*_\alpha\alpha}$ contains a
geodesic, which arises as a projection of an orbit from either a
{\it trivial} torus or a critical layer. I.e. the real homology
group $H_1(M^n,\R)$ is covered by a finite number (by clauses (i)
and (iii) of Definition 1) of corresponding homology groups of
either trivial tori or strata, forming the critical layers. That is,
if Lemma 6 was not true, it would be possible to cover the group
$H_1(M^n,\R)$ of rank $n$ by a finite number of groups of lower
dimensions, which is a contradiction.

\medskip
\noindent {\bf Proof of Lemma 7.} Consider two  Liouville tori $T_1$
and $T_2$ from the same $n$-dimensional cell of the complex $C$. As
the cell contains one non-trivial Liouville torus, all the tori in
the cell are non-trivial, by the clause (iii) of Definition 1. Let
us assume, the tori $T_1$ and $T_2$ carry the inverse images (under
the natural projection) $\gamma_1$ and $\gamma_2$ of geodesics of
class $p\alpha$ and $q\alpha$, belonging to the sets $Z_{p\alpha}$
and $Z_{q\alpha}$ respectively, with $p\neq q$. Therefore, the
(constant) frequency vectors for the phase trajectories on the tori
$T_1$ and $T_2$ should be collinear. Indeed, assuming otherwise
implies a contradiction, as $q\gamma_1(p\gamma_2)^{-1}$ is homotopic
to zero, but the tori $T_1$ and $T_2$ are non-trivial. Therefore
$\gamma_1$ and $\gamma_2$ are homotopically identical, up to
multiplicity.

\vspace{.2in} \medskip \noindent{\bf Acknowledgement:} Research
supported by EPSRC grant GR/S13682/01.

\bibliographystyle{unsrt}

\begin{thebibliography}{99}

\bibitem{A1} V.I. Arnold. Mathematical methods of classical mechanics.
Translated from the Russian by K. Vogtmann and A. Weinstein. Second
edition. Graduate Texts in Mathematics, {\bf 60}. {\em
Springer-Verlag, New York}, 1989. xvi+508 pp.

\bibitem{BT} A. Bolsinov, I.A. Taimanov. Integrable geodesic flows with positive
topological entropy.  {\em Invent. Math.}  {\bf 140} (2000), no. 3,
639--650.

\bibitem{Fom} A.T. Fomenko. Integrability and nonintegrability in geometry and mechanics. Translated from the Russian by M. V. Tsaplina.
 Mathematics and its Applications (Soviet Series), {\bf 31}. {\em Kluwer Academic Publishers Group, Dordrecht}, 1988. xvi+343 pp.

\bibitem{Fo} A.T. Fomenko. A Morse theory for integrable Hamiltonian systems. (Russian)  {\em Dokl. Akad.
Nauk SSSR}  {\bf 287}  (1986),  no. 5, 1071--1075.

\bibitem{Ko} V.V. Kozlov. Topological obstacles to the integrability of natural mechanical
systems. (Russian)  {\em Dokl. Akad. Nauk SSSR}  {\bf 249}  (1979), no. 6,
1299--1302.

\bibitem{KT} V.V. Kozlov, V.V. Ten. Topology of regions of the possibility of motion of integrable
systems. (Russian)  {\em Mat. Sb.}  {\bf 187}  (1996),  no. 5,
59--64; translation in  {\em Sb. Math.}  {\bf 187}  (1996),  no. 5,
679--684.

\bibitem{Ta} I. Taimanov. Topological properties of integrable geodesic flows. (Russian)
{\em Mat. Zametki}  {\bf 44}  (1988),  no. 2, 283--284.



\end{thebibliography}

\end{document}